\documentclass[english, 12pt,a4paper]{amsart}

\usepackage{amssymb,amsbsy,amsmath,amsfonts,amssymb,amscd,mathrsfs}
\usepackage{graphics,color,fancyhdr,multicol,multirow,fancybox}
\usepackage{wasysym,hyperref,lipsum}

\theoremstyle{plain}
\newtheorem{theorem}{Theorem}[section]

\newtheorem{corollary}[theorem]{Corollary}
\newtheorem{proposition}[theorem]{Proposition}
\newtheorem{conjecture}[theorem]{Conjecture}
\theoremstyle{definition}
\newtheorem{definition}[theorem]{Definition}

\newtheorem{remark}[theorem]{Remark}
\numberwithin{equation}{section}

\numberwithin{equation}{section}
\usepackage[top=1.5in, bottom=1.5in, left=1in, right=1in]{geometry}
\usepackage[hyperpageref]{backref}

\usepackage{xcolor}
\hypersetup{
	colorlinks,
	linkcolor={red!50!black},
	citecolor={blue!62!black},
	urlcolor={blue!80!black}
}

\def\HF{\mathop{\rm HF}\nolimits}

\def\Hom{\mathop{\rm Hom}\nolimits}

\def\Ann{\mathop{\rm Ann}\nolimits}

\def\CBP{\mathop{\rm CBP}\nolimits}
\def\span{\mathop{\rm span}\nolimits}

\def\epsilon{{\varepsilon}}

\def\phi{{\varphi}}
\let\Psi=\varPsi
\let\Phi=\varPhi
\let\theta=\vartheta

\newcommand{\X}{{\mathbb{X}}}
\newcommand{\bbP}{{\mathbb{P}}}

\newcommand{\Y}{{\mathbb{Y}}}

\newcommand{\Z}{{\mathbb{Z}}}
\newcommand{\N}{{\mathbb{N}}}

\begin{document}
		
\title{The Cayley-Bacharach property and the Levinson-Ullery conjecture}

\author{Ngoc Long Le}
\address{Department of Mathematics, University of Education - Hue University,
	34 Le Loi, Hue, Vietnam}
\email{lelong@hueuni.edu.vn}

	\author{Tran N.~K.~Linh}
	\address{Department of Mathematics, University of Education - Hue University,
		34 Le Loi, Hue, Vietnam}
	\email{tnkhanhlinh@hueuni.edu.vn}
	
	\date{\today}
	
\begin{abstract}
In this paper, we study the geometric configurations of a finite set of points
having the Cayley-Bacharach property in the $n$-dimensional projective space $\bbP^n$.
Our main contribution is the proof of the Levinson-Ullery conjecture 
for the previously unsolved case where $d=4$ and $r\ge 1$.
\end{abstract}
	
\maketitle

%
%
\section{Introduction}
\label{introduction}

Let $\X$ be a finite set of points in the $n$-dimensional projective space $\bbP^n$
over a field of characteristic zero.
We are interested in studying the intrinsic geometric properties of $\X$,
particularly, the Cayley-Bacharach property. This property is a geometric invariant
that has attracted considerable and ongoing research interest. 
Originating in classical projective geometry, the Cayley-Bacharach property 
has a long and rich history and remains a key concept in modern algebraic geometry.
Geometrically, the finite set of points $\X$ is said to have the 
\textit{Cayley-Bacharach property of degree $r$} (for short, \textit{$\X$ has $\CBP(r)$})
if every hypersurface of degree $r$ that contains all but one point of~$\X$
automatically contains the last point.  
The most celebrated classical result concerning this notion is 
the Cayley-Bacharach theorem, which asserts that if $\X$ is the complete intersection 
of two plane curves of degrees $d$ and $e$, respectively, then $\X$ has $\CBP(d+e-3)$
(see \cite{Bac, Cay2}).
Since then, this result has been generalized extensively to various contexts 
(see, e.g., \cite{DGO, EGH, GKR, Kre2, KR2, KR3, KLLT2020, KLR2019}).
Furthermore, the property has not only led to the formulation of 
many interesting conjectures (see \cite{EGH, LU2022} for historical background 
and open problems), but has also yielded significant applications to different branches
of mathematics. These applications include the algebraic characterizations of 
specific 0-dimensional schemes like complete intersections
and arithmetically Gorenstein schemes \cite{KLL2019, KLR2019, Long2023}, 
its use in studying the measures of irrationality for projective varieties 
\cite{BCD2014, LP1994, Pic2023}, and its crucial role in coding theory  
\cite{GLS2005, Han1994}.

In this paper we focus on investigating the geometrical configurations
of a finite set of points $\X\subseteq \bbP^n$ satisfying the Cayley-Bacharach property.
An interesting result on this direction is the following theorem 
(see \cite[Lemma 2.4]{BCD2014}).
\begin{theorem} \label{Thm-LineConf-BCD2014}
Let $\X\subseteq \bbP^n$ be a finite set of points having $\CBP(r)$. If 
$$
|\X| \le 2r+1
$$
then $\X$ lies on a line.
\end{theorem} 

Inspired by this result, in \cite{LU2022}, Levinson and Ullery 
introduced the notion of a plane configuration to describe the geometry 
of~$\X$ and proposed the following conjecture.
Here a \textit{plane configuration} $\mathcal{P}$ is a union of distinct
positive-dimensional linear spaces $\mathcal{P}=\bigcup_{i=1}^k\mathcal{P}_i$ 
in $\bbP^n$, its \textit{dimension} is
$\dim(\mathcal{P})=\sum_{i=1}^k\dim(\mathcal{P}_i)$ 
and its \textit{length} is $l(\mathcal{P})=k$.
\begin{conjecture}[\textbf{Levinson-Ullery}] \label{Conj-Levinson-Ullery}
	Let $\X\subseteq \bbP^n$ be a set of points having $\CBP(r)$. If 
	$$
	| \X| \le (d+1)r+1
	$$
	then $\X$ lies on a plane configuration $\mathcal{P}$ of dimension $d.$
\end{conjecture} 
Clearly, Conjecture~\ref{Conj-Levinson-Ullery} holds true for $d=1$
by Theorem~\ref{Thm-LineConf-BCD2014}. 
The main result established in the work of Levinson and Ullery \cite{LU2022} 
is summarized as follows (see \cite[Theorem 1.3]{LU2022}):
\begin{theorem} \label{Thm-Levinson-Ullery}
	The conjecture~\ref{Conj-Levinson-Ullery} holds in the following cases:
	\begin{enumerate}
		\item For all pairs $(d,r)$ with $r\le 2$ and $d\ge 1$.
		\item For all pairs $(d,r)$ with $d\le 3$ and $r\ge 1$. 
		\item For the pair $(d,r)=(4,3)$.
	\end{enumerate}
\end{theorem}
It is therefore natural to investigate whether Conjecture~\ref{Conj-Levinson-Ullery} 
remains true for all pairs $(d,r)$ with $d\ge 4$ and $r\ge 3$.
We address the case where $d=4$ and $r\ge 1$ and provide the following main result.
\begin{theorem}\label{ThmLL}
The Conjecture~\ref{Conj-Levinson-Ullery} holds for all pairs 
$(d,r)$ with $d=4$ and $r\ge 1$.
\end{theorem}
The paper is structured as follows. 
In Section 2, we establish the necessary background. 
We recall the definitions of the Cayley-Bacharach property, Hilbert functions, 
canonical modules, and plane configurations. 
We then present several algebraic characterizations of the Cayley-Bacharach property. Additionally, we include technical results based on the work of \cite{LU2022}.
Section 3 is devoted to the proof of the main theorem and a requisite proposition.

%
%
\bigbreak 
\section{Preliminary}\label{sec:0dimSchemes}
	
In what follows we work over a field $K$ of characteristic zero.
By $\bbP^n$ we denote the $n$-dimensional projective space over $K$.
Our main objects of study are finite sets of points $\X \subseteq \bbP^n$.
Let $I_\X$ denote the homogeneous vanishing ideal of $\X$ in $P=K[X_0,\dots,X_n]$,
and let $R = P/I_\X$ be the homogeneous coordinate ring of $\X$.
The ring $R$ is a 1-dimensional standard graded $K$-algebra.
After performing a homogeneous linear change of coordinates 
(see, e.g., \cite[Lemma~6.3.20]{KR2}), we may
assume that $x_0$, the image of $X_0$ in~$R$, is a non-zerodivisor of~$R$.

The Hilbert function of $\X$ is the function $\HF_\X: \Z \rightarrow \N,$
$i\mapsto \dim_K(R_i)$, and its first difference function 
$\Delta\HF_\X(i) := \HF_\X(i)-\HF_\X(i-1)$.
We have $\HF_\X(i)=0$ for $i<0$ and 
$$
1 = \HF_\X(0)<\HF_\X(1)<\cdots<\HF_\X(r_\X)=\HF_\X(r_\X+1)=\cdots=|\X|,
$$
where $r_\X$ is the regularity index of $\X$.

\begin{definition}
For $r\ge 0$, we say that $\X$ has the \textbf{Cayley-Bacharach property 
of degree $r$} (for short, $\X$ has $\CBP(r)$) if every hypersurface of degree $r$
which contains all but one point of $\X$ must contain all the points of $\X$.
In the case that $\X$ has $\CBP(r_\X-1)$ we also say that $\X$ is a 
\textbf{Cayley-Bacharach scheme}.
\end{definition}

Notice that $\X$ has $\CBP(r)$ if and only if for any $p\in\X$
we have $\HF_{\X\setminus\{p\}}(r)=\HF_\X(r)$.
If $\X$ has $\CBP(r)$ then it has $\CBP(r-1)$.
Moreover, the number $r_\X-1$ is the largest degree $r\ge 0$ such that 
$\X$ can have $\CBP(r)$.

\begin{remark}
Let $\Y$ be a subset of $\X$ with $|\Y| =|\X|-1$.
\begin{enumerate}
	\item[(a)] The image of the vanishing ideal $I_\Y\subseteq P$ of $\Y$ 
	in $R=P/I_\X$ is denoted by $I_{\Y/\X}$.
	Its initial degree is $\alpha_{\Y/\X} := \min\{k\in\N \mid (I_{\Y/\X})_k\ne 0\}$.
	Then $\alpha_{\Y/\X}\le r_\X$ and the Hilbert function of $\Y$ satisfies
	$$
	\HF_{\Y}(i) \;=\;
	\begin{cases}
		\HF_\X(i) & \text{for $i<\alpha_{\Y/\X}$},\\
		\HF_\X(i)-1 & \text{for $i\ge\alpha_{\Y/\X}$}.
	\end{cases}
	$$
	
	\item[(b)] A nonzero homogeneous element $f^*_\Y \in (I_{\Y/\X})_{\alpha_{\Y/\X}}$
		is called a \textbf{minimal separator} of $\Y$ in $\X$.
		In this case, $(I_{\Y/\X})_{\alpha_{\Y/\X}+i} = Kx_0^i f^*_\Y$ 
		for every $i\ge 0$.	
		
	\item[(c)] The graded $R$-module $\omega_R = \Hom_{K[x_0]}(R,K[x_0])(-1)$
	is called the \textbf{canonical module} of $R$. Its $R$-module structure is defined by
	$(f\cdot \varphi)(g)=\varphi(fg)$ for all $f,g\in R$ and $\varphi\in \omega_R$.
	It is also a finitely generated graded $R$-module and for every $i\in\Z$ we have
	$$
	\HF_{\omega_R}(i) =|\X| -\HF_\X(-i)
	$$
	(see \cite[Lemma~1.3]{Kre2}).
\end{enumerate}
\end{remark}

From \cite[Proposition~2.1]{GKR} and \cite[Proposition~4.3]{KLLT2020}
we have the following characterizations of the Cayley-Bacharach property of $\X$.

\begin{proposition}\label{CB-Characterizations}
The following conditions are equivalent for every $r\ge 1$:
\begin{enumerate}
	\item[(a)] $\X$ has $\CBP(r)$.
	
	\item[(b)] Every subset $\Y$ of $\X$ with $|\Y|=|\X|-1$ 
	satisfies $\alpha_{\Y/\X}\ge r+1$.
	
	\item[(c)] No element of $(I_{\Y/\X})_{r_\X}\setminus\{0\}$ is divisible by $x_0^{r_\X-r}$
	
	\item[(d)] There exists an element $\varphi\in (\omega_R)_{-r}$ such that
	$\Ann_R(\varphi)=0$.
\end{enumerate}
\end{proposition}

\begin{corollary}\label{Cor-LowerBoundCBP}
Let $\X\subseteq \bbP^n$ be a finite set of points having $\CBP(r)$.
\begin{enumerate}
	\item[(a)] $\HF_\X(i)+\HF_\X(r-i)\le |X|$ for all $i\in\Z$.
	
	\item[(b)] We have $|\X| \ge r+2$.
\end{enumerate}
\end{corollary}
\begin{proof}
By Proposition~\ref{CB-Characterizations}(d), there 
exists $\varphi\in (\omega_R)_{-r}$ such that $\Ann_R(\varphi)=0$.
The multiplication map $R\rightarrow R\cdot\varphi \subseteq \omega_R(-r), f\mapsto f\cdot\varphi$ 
is injective, and hence 
$$
\HF_\X(i) \le \HF_{\omega_R(-r)}(i)  = \HF_{\omega_R}(i-r) 
=|\X| -\HF_\X(r-i)
$$
for all $i\in\Z$. This completes the proof of (a).
Claim (b) follows from (a) or \cite[Lemma~2.4]{BCD2014}.
\end{proof}

Now we recall from \cite[Definition~2.1]{LU2022} the following notion of plane configurations.

\begin{definition}
A \textbf{plane configuration} $\mathcal{P}$ is a union of distinct positive-dimensional 
linear spaces $\mathcal{P}=\bigcup_{i=1}^k\mathcal{P}_i$ in~$\bbP^n$.
The number $\dim(\mathcal{P}) :=\sum_{i=1}^k\dim(\mathcal{P}_i)$ is called 
the \textbf{dimension} of $\mathcal{P}$ and $\ell(\mathcal{P}) := k$ 
is called the \textbf{length} of $\mathcal{P}$.
\end{definition}

Note that a 0-dimensional plane configuration is empty, a 1-dimensional plane
configuration is a single line, and a 2-dimensional plane configuration is either two 
lines or one 2-plane. Here a $d$-plane refers to a $d$-dimensional linear space.
Also, the linear span of a plane configuration 
$\mathcal{P}=\bigcup_{i=1}^k\mathcal{P}_i$ is denoted by
$\span(\mathcal{P})$ or $\span(\mathcal{P}_i\mid i=1,\dots,k)$.  

\begin{definition}
Let $\mathcal{P} =\bigcup_{i=1}^k\mathcal{P}_i$ be a plane configuration in $\bbP^n$.
\begin{enumerate}
	\item[(a)] $\mathcal{P}$ is said to be \textbf{skew} 
	if $\mathcal{P}_i\cap\mathcal{P}_j =\emptyset$ for all $i\ne j$.
	
	\item[(b)] $\mathcal{P}$ is said to be \textbf{split} if
	$\mathcal{P}_i\cap \span(\mathcal{P}_j \mid j\ne i) =\emptyset$
	for every $i=1,\dots,k$.
\end{enumerate}
\end{definition}

Clearly, if $\mathcal{P}$ is split then it is skew, and the converse holds true 
if $\ell(\mathcal{P})=2$. Also, $\mathcal{P}$ is split if and only if 
$\dim(\span(\mathcal{P})) =\dim(\mathcal{P})+\ell(\mathcal{P})-1$. 

The following proposition follows from \cite[Proposition~2.5]{LU2022}.

\begin{proposition}\label{Prop-CBP-Complement}
Let $\X \subseteq \bbP^n$ be a finite set of points, and
let $\mathcal{P}$ be a plane configuration of length $\ell(\mathcal{P})=k$.
If $\X$ has $\CBP(r)$ then $\X \setminus \mathcal{P}$ has $\CBP(r-k)$.
\end{proposition}
	
The results presented in \cite[Propositions~4.2, 4.7 and Lemma~4.5(1)]{LU2022} 
lead to the next proposition.

\begin{proposition} \label{CBP-PlaneConfiguration} 
Let $\X \subseteq \bbP^n$ be a finite set of points, 
and let $\mathcal{P}=\bigcup_{i=1}^k\mathcal{P}_i$ be a plane configuration of length
$\ell(\mathcal{P})=k$ such that $\X_{\mathcal{P}_i} := \X\cap \mathcal{P}_i\ne \emptyset$ 
for all $i=1,\dots,k$.
\begin{enumerate}
	\item[(a)]  
	If $\mathcal{P}$ is split and $\X\subseteq \mathcal{P}$,
	then $\X$ has $\CBP(r)$ if and only if, for
	each $i$, $\X\cap \mathcal{P}_i$ has $\CBP(r)$.
	
	\item[(b)] 
	If $\X$ has $\CBP(r)$ and $\mathcal{P}$ is skew, then 
	one of the following holds:
	\begin{enumerate}
		\item[(i)] Each plane $\mathcal{P}_i$ contains at least 
		$\max\{k, r + 2)$ points of $\X$, or
		
		\item[(ii)] Some plane $\mathcal{P}_i$ contains fewer than $k$ points, 
		and also $k \ge r+2.$
	\end{enumerate}
	
	\item[(c)] 
	Suppose that $\X$ has $\CBP(r)$, $k=2$, and 
	$\mathcal{P}_1$ and $\mathcal{P}_2$ meet at a point $p\in\bbP^n$. Then at least 
	one of $\X_{\mathcal{P}_1}$ and $\X_{\mathcal{P}_1}\cup\{p\}$ has $\CBP(r).$
\end{enumerate}
\end{proposition}

%
%
\bigbreak 
\section{Proof of main theorem}

As in the previous section, let $r\ge 1$ and let $\X\subseteq \bbP^n$ be a finite set 
of points having $\CBP(r)$. This section is devoted to the proof of Theorem~\ref{ThmLL}, 
which states that if $|\X| \le 5r+1$ then $\X$ lies on a plane 
configuration $\mathcal{P}$ of dimension $4$.  
This proof will utilize the following proposition.

\begin{proposition}\label{Prop-InductiveTechnique}
	Let $\X\subseteq\bbP^n$ be a finite set of points having $\CBP(r)$ and 
	$|\X|\le (d+1)r+1$, let $A$ be a $k$-plane with $k\le d$, 
	and let $\X_A := \X \cap A$ and $\X_B := \X\setminus \X_A$. 
	\begin{enumerate}
		\item[(a)] \label{PropLB}
		If Conjecture \ref{Conj-Levinson-Ullery} holds for $(d-1,r)$
		and $\X$ does not lie on a plane configuration of dimension $d-1$, 
		then $|\X| \ge dr+2$.
		
		\item[(b)] \label{PlaneDecomposition}
		Assume that Conjecture \ref{Conj-Levinson-Ullery} holds 
		for every pair $(i,j)$ with $i\le d$ and $j<r$.
		\begin{enumerate}
			\item[(i)] If $|\X_{A}| \ge d+1$ then $\X_B$ lies on 
			a plane configuration of dimension $d$.
			
			\item[(ii)] If $\X_A$ does not lie on a plane configuration 
			of dimension $(d-1)$ and $\X_B$ is contained in an $\ell$-plane 
			with $\ell\le d$ and $d\le r$, then $\X = \X_A$.
		\end{enumerate}
	\end{enumerate}
\end{proposition}
\begin{proof}
	(a)\quad Note that $\X$ has $\CBP(r)$.
	If $|\X| \le dr+1$, then the validity of Conjecture~\ref{Conj-Levinson-Ullery} 
	for the pair $(d-1,r)$ guarantees that $\X$ lies on a plane configuration of dimension $d-1$.
	
	\medskip\noindent(b)\quad
	If $\X_B=\emptyset$ then we are done.
	Suppose that $\X_B\ne \emptyset$. By Proposition~\ref{Prop-CBP-Complement}, 
	$\X_B$ has $\CBP(r-1)$.
	Since $|\X_{A}|\ge d+1$, we have 
	$$
	1\le |\X_B| \le (d+1)r+1 -(d+1) =(d+1)(r-1)+1.
	$$
	Because Conjecture~\ref{PlaneDecomposition} holds true for $(d,r-1)$, 
	the set $\X_B$ must lie on a plane configuration of dimension~$d$.
	Thus claim (i) follows.
	
	Next, we prove (ii). 
	Observe that $A$ is a $d$-plane and $|\X_A|\ge d+1$. 
	Suppose for a contradiction that $\X_B\ne \emptyset$.
	Let $B$ be an $\ell$-plane containing $\X_B$.
	By Proposition~\ref{Prop-CBP-Complement}, both $\X_A$ and $\X_B$ have $\CBP(r-1)$.
	Since $\X_A$ does not lie on a plane configuration of dimension $(d-1)$
	and Conjecture~\ref{Conj-Levinson-Ullery} holds for $(i,j)$ with $i\le d$ 
	and $j\le r-1$, claim (a) yields  
	$$
	|\X_A| \;\ge\; d(r-1)+2.
	$$
	It follows that 
	$$
	|\X_B|\le ((d+1)r+1)-(d(r-1)+2) = r+d-1 \le 2(r-1)+1,
	$$
	and consequently $\X_B$ lies on a line $L$ by Theorem~\ref{Thm-LineConf-BCD2014}. 
	Now we distinguish the following two cases.
	\begin{enumerate}
		\item[(i)] \textit{$A$ and $L$ are split.} In this case, $\X_A$ and $\X_B$
		have $\CBP(r)$ by Proposition~\ref{CBP-PlaneConfiguration}(a).
		Due to Corollary~\ref{Cor-LowerBoundCBP} and (a), 
		the lower bounds for the cardinalities of $\X_A$ and $\X_B$ are given by 
		$$
		|\X_A| \ge dr+2\quad\text{and}\quad
		|\X_B| \ge r+2.
		$$
		This implies $(d+1)r+4 \le |\X_A|+|\X_B| = |\X|\le (d+1)r+1$, 
		which is impossible.
		
		\item[(ii)] \textit{$A$ and $L$ intersect at a point, say $p$.} 
		By Lemma~\ref{CBP-PlaneConfiguration}(c), at least one of 
		$\X_A$ and $\X_A\cup \{p\}$ has $\CBP(r).$
		If $\X_A$ has $\CBP(r)$, then claim (a) yields that
		$|\X_A| \ge dr+2$. When $\X_A\cup \{p\}$ has $\CBP(r),$ then 
		$|\X_A \cup \{p\}| \ge dr+2$ by (a). 
		In these both cases, we always have $|\X_A| \ge dr+1$. 
		On the other hand, the set $\X_B$ has $\CBP(r-1)$, and Corollary~\ref{Cor-LowerBoundCBP}
		shows that $|\X_B| \ge r+1$.  
		Hence, we obtain $|\X| = |\X_A|+|\X_B| \ge (dr+1)+(r+1) \ge (d+1)r+2$, 
		a contradiction.
	\end{enumerate}	
	Altogether, we find $\X_B=\emptyset$ and $\X = \X_A$.
\end{proof}

Now we are ready to give a proof of our main result.

\begin{proof}[Proof of Theorem~\ref{ThmLL}]
According to Theorem~\ref{Thm-Levinson-Ullery}, it suffices to consider the case $r\ge 4$
and use induction on $r$.
Let $\X\subseteq \bbP^n$ be a finite set of points having $\CBP(r)$ and
$|\X| \le 5r+1$. Further, let $A$ be a 4-plane containing the maximum number 
of points of $\X$, let $\X_A :=\X\cap A$, and let 
$\X_B := \X\setminus \X_A$.
When $\X_B=\emptyset$, the theorem holds true.
So, we consider the case $\X_B\ne \emptyset$.
In this case $|\X_A| \ge 5$, and Proposition~\ref{Prop-InductiveTechnique}(b.i) 
yields that the set $\X_B$ lies on a plane configuration of dimension 4.
Also, the set $\X_B$ has $\CBP(r-1)$ by Proposition~\ref{Prop-CBP-Complement}. 
In the following we consider different
cases depending on the geometric configurations of $\X_B$.
\begin{enumerate}
	\item[\textit{Case 1:}] \textit{$\X_B$ lies on a 4-plane.} \\
	If $\X_A$ does not lie on a plane configuration of dimension 3, then,
	by the inductive hypothesis, Proposition~\ref{Prop-InductiveTechnique}(b.ii) 
	shows that $\X = \X_A$, and subsequently $\X$ lies on a plane configuration 
	of dimension 4, as wanted.
	By the choice of $A$, $\X_A$ lies on a 4-plane, but not a 3-plane.
	Thus, it suffices to look at the case that $\X_A$ lies on a plane configuration
	of dimension 3, i.e., it lies on either a skew plane configuration
	of three lines or a split plane configuration of a 2-plane and a line. 
    \begin{enumerate}
	\item[(1.1)] \textit{$\X_{A}$ lies on a skew plane configuration 
	of three lines $L_1, L_2, L_3$.} \\
	Note that $\X_A$ has $\CBP(r-1)$ by Proposition~\ref{Prop-CBP-Complement} 
	and $r\ge 4$. By Proposition~\ref{CBP-PlaneConfiguration}(b),  each line $L_i$ contains at least $r+1$ points of
	$\X_A$ for $i=1,2,3$. Consequently, we have $|\X_A| \ge 3(r+1)$.
	It follows that 
	$$
	|\X_B| \le 5r+1 - 3(r+1) = 2r-2< 2(r-1)+1.
	$$
	According to Theorem \ref{Thm-LineConf-BCD2014}, the set $\X_{B}$ lies on a line,
	and therefore, the set $\X$ lies on a plane configuration of dimension 4,
	as desired. 

	\item[(1.2)] \textit{$\X_{A}$ lies on a split plane configuration 
		of a 2-plane $H$ and a line $L$.}\\
	Since $\X_{A}$ has  $\CBP(r-1)$, Proposition~\ref{CBP-PlaneConfiguration}(a) implies
	that two sets of points $\X_A\cap H$ and $\X_A\cap L$ also have $\CBP(r-1)$.
	Moreover, $\X_A\cap H$ does not lie on a line by the choice of $A$.
	By Corollary~\ref{Cor-LowerBoundCBP} and
	Proposition~\ref{Prop-InductiveTechnique}(a), we get 
	$$
	|\X_{A}\cap L|\ge r+1\quad\text{and}\quad
	|\X_{A}\cap H|\ge 2(r-1)+2=2r.
	$$
	This implies 
	$$
	|\X_{B}|\le 5r+1-(r+1)-2r\le 2(r-1)+2.
	$$ 
		\begin{enumerate}
		\item[(i)] If $|\X_{B}|\le 2(r-1)+1$ then $\X_{B}$ must lies on a line 
		by Theorem~\ref{Thm-LineConf-BCD2014}, and hence the set $\X$ is contained 
		in a plane configuration of dimension 4, as wanted to show.
		
		\item[(ii)] If $|\X_{B}|=2r \le 3(r-1)+1$, then $\X_{B}$ lies on 
		a plane configuration of dimension 2 by the inductive hypothesis 
		and $|\X_A|= 3r+1$.
		When $\X_B$ lies on a line, as above, the set $\X$ lies on 
		a plane configuration of dimension 4 and we are done.
		If $\X_B$ lies on a union of two split lines, then each of these lines 
		contains at least $r+1$ points of $\X_B$ 
		by Proposition~\ref{CBP-PlaneConfiguration}(a) and Corollary~\ref{Cor-LowerBoundCBP}, and so $|\X_B|\ge 2r+2$,
		which is impossible. Here we remark that a skew plane configuration
		of two lines agrees with a split plane configuration of two lines.
		Now we consider the last subcase that 
		the set $\X_B$ lies on a $2$-plane~$H'$. 
		If $H'\cap H\ne \emptyset$ or $H'\cap L\ne \emptyset$ then
		the pair $(H,H')$ defines a 4-plane containing at least $4r$ points of $\X$
		or the pair $(L,H')$ and a point in $\X_A\cap H$ define
		a 4-plane containing at least $3r+2$ points of $\X$,
		which contradicts to the choice of the 4-plane $A$.
		Hence we must have $H'\cap H = H'\cap L =\emptyset$. 
		In particular, $\X$ is contained in a skew plane configuration 
		$\mathcal{P}=H\cup L\cup H'$ of length $\ell(\mathcal{P})=3 < r$.
		Since $\X$ has $\CBP(r)$, an application 
		of Proposition~\ref{CBP-PlaneConfiguration}(b)
		yields that the line $L$ contains at least $r+2$ points of $\X$.  
		It follows that 
		$$
		|\X| = |\X\cap H| + |\X\cap L| +|\X\cap H'| 
		\ge 2r+(r+2)+2r = 5r+2.
		$$
		However, this contradicts the assumption that $|\X|\le 5r+1$.
		\end{enumerate}
	\end{enumerate}
\end{enumerate}
Next, we assume that no 4-plane contains $\X_B$
and we consider the following remaining cases.
\begin{enumerate}
\item[\textit{Case 2}:] \label{Case2}	
\textit{$\X_B$ is contained in the union of a 3-plane $H$ and a line $L$.} \\
Since $\X_B$ is not contained in a $4$-plane, then the 3-plane $H$ and 
the line $L$ must be split. We define two subsets of $\X_B$ by 
$\X_H:=\X_B\cap H$ and $\X_L:= \X_B\cap L$.
Remark that $\X_B$ has $\CBP(r-1)$. By Proposition~\ref{CBP-PlaneConfiguration}(a), 
two subsets $\X_H$ and $\X_L$ also have $\CBP(r-1)$. 
By Corollary~\ref{Cor-LowerBoundCBP}, we have $|\X_L| \ge r+1$.
Observe that $\X_H$ may lie on a plane configuration of dimension 2,
but it does not lie on a line.  
By Proposition~\ref{Prop-InductiveTechnique}(a), the cardinality of $\X_H$ is bounded below as
$$
| \X_H| \ge 2(r-1)+2 = 2r.
$$
Moreover, since the 4-plane $A$ contains the maximum number of points in $\X$, 
we must have 
$$
|\X_A| \ge |\X_H| +1 \ge 2r + 1.
$$ 
Thus, combining the lower bounds for the disjoint sets 
$\X_A, \X_H, \X_L$, we find
$$
|\X| = |\X_A| + |\X_H| + |\X_L| \ge 
(2r+1)+2r + r+1 = 5r+2,
$$ 
which contradicts the hypothesis that $|\X|\le 5r+1$.
\item[\textit{Case 3:}] \label{Case3} 
\textit{$\X_B$ is contained in the union of two 2-planes $K_1$ and $K_2$
such that $\X_{K_i} := \X_B\cap K_i$ does not lie on a line for $i=1,2$.} \\
Clearly, the two planes $K_1$ and $K_2$ are split, as $\X_B$ is not contained 
in a $4$-plane.  Also, we see that $\X_B=\X_{K_1}\cup \X_{K_2}$
and $\X_B$ has $\CBP(r-1)$.
It follows from Proposition~\ref{CBP-PlaneConfiguration}(a) that 
both $\X_{K_1}$ and $\X_{K_2}$ have $\CBP(r-1)$.
Thank to Proposition~\ref{Prop-InductiveTechnique}(a), the cardinality of $\X_B$ satisfies
$$
|\X_B| = |\X_{K_1}| +|\X_{K_2}| \ge 2r+2r = 4r.
$$
By the choice of the 4-plane $A$, we must have
$$
|\X_A| \ge |\X_{K_1}|+2 \ge 2r+2.
$$
Hence we get $|\X| = |\X_A| + |\X_B| \ge 4r+2r+2 =6r+2 >5r+1$,
which is impossible. 
\item[\textit{Case 4:}] \label{Case4} 
\textit{$\X_B$ is contained in the union of a $2$-plane $K$ and two lines 
$L_1$ and $L_2$ such that $\X_K := \X_B\cap K$ does not lie on a line. }\\
In this case $L_1$ and $L_2$ must be skew, since otherwise we would go back
to Case~3 by defining the 2-plane $K' =\span(L_1, L_2)$.
If $K\cap L_i\ne\emptyset$ for $i\in\{1,2\}$ then we go back to Case~2.
Thus $K\cup L_1\cup L_2$ is a skew plane configuration containing $\X_B$
and $\X_B= \X_K\cup \X_{L_1}\cup \X_{L_2}$, where
$\X_{L_1} := \X_B\cap L_1$ and $\X_{L_2} :=\X_B\cap L_2$.
The linear space $H :=\span(L_1, L_2)$ is a $3$-plane. Put $\X_H:=\X_B\cap H$.
We distinguish the following subcases.
\begin{enumerate}
\item[(4.1)] $H\cap K=\emptyset$:\quad  By Proposition~\ref{CBP-PlaneConfiguration}(a),
the sets $\X_H$ and $\X_K$ have $\CBP(r-1)$, since $\X_B$ has $\CBP(r-1)$.
From $\X_H = \X_{L_1}\cup \X_{L_2}$ and $L_1\cap L_2=\emptyset$,
again Proposition~\ref{CBP-PlaneConfiguration}(a) yields that  
$\X_{L_1}$ and $\X_{L_2}$ have $\CBP(r-1)$. 
By Corollary~\ref{Cor-LowerBoundCBP} and Proposition~\ref{Prop-InductiveTechnique}(a), 
we have 
$$
|\X_{K}|\ge 2(r-1)+2,\quad |\X_{L_1}|\ge r+1,\quad  |\X_{L_2}|\ge r+1,
$$
and consequently 
$$
|\X_B| = |\X_{K}|+ |\X_{L_1}|+ |\X_{L_2}| \ge 2r+2(r+1)=4r+2.
$$
Furthermore, by the choice of the 4-plane $A$, we find
$$
|\X_A| \ge |\X_B\cap H|+1 = |\X_{L_1}|+ |\X_{L_2}| +1\ge 2r+3.
$$ 
Thus, we obtain $6r+5\le |\X_A| +|\X_B| = |\X|\le 5r+1$, which
is impossible.
\item[(4.2)] $H\cap K\ne \emptyset$:\quad 
Because the set $\X_B$ is not contained in a $4$-plane, the 3-plane $H$ 
intersects the 2-plane $K$ at a point, say $p$. 
By Proposition~\ref{CBP-PlaneConfiguration}(c), $\X_{K}$ or $\X_{K}\cup \{p\}$ 
has $\CBP(r-1)$. If $\X_{K}$ has $\CBP(r-1)$ then by
Proposition~\ref{Prop-InductiveTechnique}(a) we have
$$
| \X_K| \ge 2(r-1)+2 =2r.
$$
When $\X_{K}\cup \{p\}$ has $\CBP(r-1)$, we have
$$
| \X_K| \ge 2(r-1)+1 = 2r-1.
$$
This yields $|\X_K|\ge 2r-1$. On the other hand, we have
$\X_{L_1}\cup \X_{L_2} = \X_B\setminus K$, and so 
this union has $\CBP(r-2)$. Especially, $\X_{L_1}$ and $\X_{L_2}$
also have $\CBP(r-2)$. By Corollary~\ref{Cor-LowerBoundCBP}, we find
$|\X_{L_i}| \ge (r-2)+2=r$ for $i=1,2$. So, we get
$$
|\X_H| = | \X_{L_1}| +| \X_{L_2}| +1 \ge 2r +1.
$$
This implies  $|\X_A| \ge |\X_H|+1\ge  2r+2$,
and consequently 
$$
6r+1=(2r+2)+(4r-1)\le |\X_A| +|\X_B| = |\X|\le 5r+1,
$$
which is impossible. 
\end{enumerate}
\item[\textit{Case 5:}] 
\textit{$\X_B$ is contained in the union of $s$ lines $L_1, \dots, L_s$, 
where $s=3,4$.}\\
When $s=3$, three lines $L_1, L_2, L_3$ must be split, because $\X_B$
is not contained in a 4-plane. Then $H:=\span(L_1,L_2)$ is a 3-plane, and 
$H$ and $L_3$ are split. In this situation, we go back to Case~2.
Now suppose that $s=4$.
If any two lines of these lines is contained in a 2-plane then we turn back to 
either Case 3 or Case 4. This enables us to assume that those lines are pairwise skew.
In particular, $\X_B$ does not lie on a plane configuration of dimension 3.
For $i\in\{1,\dots,4\}$, put $\X_{L_i} := \X_B\cap L_i$ and 
$H_{ij} := \span(L_i,L_j)$ for $i\ne j$.  
Note that $\X_{L_i} \ne \emptyset$ for $i=1,\dots,4$.
Consider the following subcases.
\begin{enumerate}
	\item[(5.1)]  \textit{There are three of four lines that are split.}
	Without loss of generality, we may assume that $L_1, L_2$ and $L_3$ are split. 
	Then $\X_{L_1}\cup\X_{L_2}\cup\X_{L_3} = \X_B\setminus L_4$
	and this union has $\CBP(r-2)$. By Proposition~\ref{CBP-PlaneConfiguration}(a), 
	the sets of points $\X_{L_1}, \X_{L_2}, \X_{L_3}$ have $\CBP(r-2)$. 
	It follows that $|\X_{L_i}| \ge r$ for $i=1,2,3$ 
	by Corollary~\ref{Cor-LowerBoundCBP}.
	So, we have 
	$$
	|\X_{B}| = |\X_{L_1}| +|\X_{L_2}| +|\X_{L_3}| +|\X_{L_4}|
	\ge 3r +1.
	$$
	Moreover, since $H_{12}$ is a 3-plane containing $L_1$ and $L_2$, we also have
	$$
	|\X_{A}| \ge |\X_{L_1}|+|\X_{L_2}|+1 \ge 2r+1.
	$$ 
	Consequently, in this subcase, we get 
	$$
	|\X| = |\X_{A}| + |\X_{B}| \ge 5r+2,
	$$ 
	which contradicts the hypothesis $|\X|\le 5r+1$.
	\item[(5.2)] \textit{All three of four lines are not split.}
	Without lost of generality, we may assume that 
	three lines $L_1, L_2, L_3$ satisfy the condition that
	their union $\X_{L_1}\cup \X_{L_2}\cup\X_{L_3}$ 
	contains at least $\frac{3}{4}$ points of $\X_{B}$.
	Also, we have $L_3\cap H_{12} \ne \emptyset$.
	\begin{itemize}
		\item If $L_3 \subseteq H_{12}$ then we come back to Case 2.
		
		\item If $L_3$ intersects $H_{12}$ at a point, then 
		$\span(H_{12},L_3) = \span(L_1,L_2,L_3)$ is a 4-plane.
		By the choice of the 4-plane $A$, we have 
		$$
		|\X_{A}|\ge |\X_{L_1}| + |\X_{L_2}| + |\X_{L_3}| 
		\ge \frac{3}{4}|\X_{B}|.
		$$ 
		This implies
		$$
		5r+1\ge |\X| = |\X_{A}|+|\X_{B}|\ge \frac{7}{4}|\X_{B}|.
		$$
		Thus, for $r\ge 4$, we get 
		$$
		|\X_{B}|\le \frac{20r+4}{7}\le 4(r-1)+1.
		$$
		By Theorem~\ref{Thm-Levinson-Ullery}, $\X_B$ lies on a plane configuration
		of dimension 3, and this is impossible.
	\end{itemize}
\end{enumerate}  
\end{enumerate} 
Altogether, the proof of the theorem is complete. 
\end{proof}

%
%
\goodbreak
\section*{Acknowledgements}

This work is supported by the \textit{Vietnam Ministry of Education and Training} 
under the grant number B2026-DHH-01. 

\goodbreak

\address

\end{document}